\documentclass[12pt, leqno]{amsart}

\usepackage[OT2,T1]{fontenc}
\usepackage{indentfirst}
\usepackage{amstext}
\usepackage{amsthm}
\usepackage{amsopn}
\usepackage{amsfonts}
\usepackage{amsmath}
\usepackage{latexsym}
\usepackage[francais,english]{babel}
\usepackage{amscd}
\usepackage{amssymb}
\usepackage{amsmath}
\usepackage[all,cmtip]{xy}
\usepackage{leftidx}
\usepackage{graphicx}
\usepackage{etoolbox}
\patchcmd{\section}{\normalfont\scshape\centering}{\normalfont\bfseries}{}{}
\patchcmd{\subsection}{-.5em}{.5em}{}{}

\newtheorem{theo}{{Theorem}}[section]
\newtheorem{coro}[theo]{{Corollary}}

\newtheorem{prop}[theo]{Proposition}

\theoremstyle{definition}

\newtheorem{defn}[theo]{Definition}

\numberwithin{equation}{section}

\newtheorem{notation}[theo]{Notation}

\DeclareFontEncoding{OT2}{}{} 
  \newcommand{\textcyr}[1]{%
    {\fontencoding{OT2}\fontfamily{wncyr}\fontseries{m}\fontshape{n}%
     \selectfont #1}}
\newcommand{\sha}{{\mbox{\textcyr{Sh}}}}

\begin{document}
\tolerance 400 \pretolerance 200 \selectlanguage{english}

\title{Embeddings of maximal tori in classical groups and Hasse principles}
\author{Eva  Bayer-Fluckiger}
\date{\today}
\maketitle


\medskip


\small{} \normalsize

\medskip

\selectlanguage{english}
\section{Introduction}

The aim of this paper is to revisit a topic investigated in \cite {PR}, \cite {GR}, \cite {Lee}, \cite {BLP3}, \cite{BLP2},  \cite {BLP1}, namely the question of embeddings
of maximal tori in classical groups. As in the above references, this is reformulated in terms of embeddings of 
commutative \'etale algebras with involution in central simple algebras with involution. 

\medskip
If the ground field is a global field of characteristic $\not = 2$, a necessary and sufficient condition for the local-global
principle to hold is given in \cite {BLP1}, Theorem 4.6.1. This result is formulated in terms of an obstruction group. The
first aim of the present paper is to give a simpler version of Theorem 4.6.1 of \cite {BLP1}, in terms of a different 
description of the obstruction group, see 
Theorem \ref {main}.

\medskip
We show that if the commutative \'etale algebra is a product of pairwise
independent fields,
then the obstruction group is trivial, and hence the local-global principle holds
(see Corollary \ref {pairwise}).

\bigskip
\section {Definitions, notation and basic facts} \label{definitions}

\medskip
Let $L$ be a field with ${\rm char}(L) \not = 2$, and let $K$ be a subfield of $L$ such that either $K = L$, or
$L$ is a quadratic extension of $K$. 

\medskip
{\bf \'Etale algebras with involution}

\medskip
Let $E$ be a commutative \'etale algebra of finite rank over $L$, and let $\sigma: E \to E$ be a $K$--linear  involution. Set  $F = \{e \in E | \sigma(e) = e \}$, and $n = {\rm dim}_L(E)$. Assume that
if $L = K$, then we have ${\rm dim}_K(F) = [{{n+1}\over 2}] $. Note that if $L \not = K$, then  ${\rm dim}_K(F) = n$ (cf. \cite{PR},
Proposition 2.1.).

\medskip

{\bf Central simple algebras with involution}

\medskip

Let  $A$ be a central simple algebra  over $L$. Let $\tau$ be an involution of $A$, and assume that  $K$ is the fixed field of $\tau$ in $L$. Recall that
$\tau$ is said to be of the {\it first kind} if $K = L$ and of the {\it second kind} if $K \not = L$; in this case, $L$ is a quadratic extension of $K$. 
After extension to a splitting field of $A$, any involution of the first kind is induced by a symmetric or by a skew-symmetric form. We say that the involution is of the {\it orthogonal type} in the first case, and of the {\it symplectic type} in the second case. An involution of the second kind is also called {\it unitary involution}.

\medskip
{\bf Embeddings of algebras with involution}

\medskip
Let $(E,\sigma)$ and $(A,\tau)$ be as above, with $n = {\rm dim}_L(E)$ and ${\rm dim}_L(A) = n^2$; assume 
moreover that  $\sigma | L = \tau | L$.

\medskip

An {\it embedding} of $(E,\sigma)$ in $(A,\tau)$ is by definition an injective homomorphism
$f : E \to A$ such that $\tau (f(e)) = f(\sigma(e))$ for all $e \in E$. It is well--known that embeddings
of maximal tori into classical groups can be described in terms of embeddings of \'etale algebras with involution 
into central simple algebras with involution satisfying the above dimension hypothesis (see for
instance \cite{PR}, Proposition 2.3). 

\medskip
Let $\epsilon : E \to A$ be an $L$--embedding which may not respect the given involutions. There exists an
involution $\theta$ of $A$ of the same type (orthogonal, symplectic or unitary) as $\tau$ such that
 $\epsilon (\sigma(e)) = \theta (\epsilon (e)) \ {\rm for \ all} \   e \in E,$
in other words $\epsilon : (E,\sigma) \to (A,\theta)$ is an $L$---embedding of algebras with involution
(see \cite{K}, \S2.5. or \cite{PR}, Proposition 3.1).

\medskip
For all $a \in F^{\times}$, let $\theta_a : A \to A$ be the involution given by $\theta_a =  \theta  \circ {\rm Int}(\epsilon (a))$. Note that
$\epsilon : (E,\sigma) \to (A,\theta_a)$ is an embedding of algebras with involution.

\begin{prop}\label{embedding}The following conditions are equivalent :

\smallskip 
{\rm (a)} There exists an $L$--embedding  $\iota : (E,\sigma) \to (A,\tau)$ of algebras with involution.

\smallskip 
{\rm (b)} There exists an $a \in F^{\times}$ such that $(A,\theta_a) \simeq (A,\tau)$ as algebras with involution.

\end{prop}

\noindent
{\bf Proof.} See \cite{PR} Theorem 3.2.

\medskip
{\bf Oriented embeddings}

\medskip An orthogonal involution $(A,\tau)$ is said to be {\it split} if $A$ is a matrix algebra. 
Note that if $n$ is odd, then every orthogonal involution is split. In the case of
a nonsplit orthogonal involution, we need an additional notion, called {\it orientation} (see \cite {BLP1}, \S 2).

\medskip
Assume that  $(A,\tau)$ is of orthogonal type and that $n$ is even. We denote by by $C(A,\tau)$ its {\it Clifford algebra} (cf. \cite{KMRT} Chap II.  (8.7)), and by $Z(A,\tau)$ the center of the algebra $C(A,\tau)$. Then 
$Z(A,\tau)$ is a quadratic \'etale algebra over $K$. 

\medskip
Let $\Delta (E)$ be the discriminant algebra of $E$ (cf. \cite {KMRT} Chapter V, \S 18,
p. 290). An isomorphism of $K$--algebras $$\Delta(E) \to Z(A,\tau)$$ will be called an {\it orientation} (see
\cite {BLP1}, \S 2). 

\medskip
 Let $u : \Delta(E) \to Z(A,\theta)$ be the orientation of $(A,\theta)$ constructed in  \cite {BLP1}, 2.3, and for all $a \in F^{\times}$ let $u_a : \Delta(E) \to Z(A,\theta_a)$ be as in \cite {BLP1}, 2.5. 
Recall the notion of  oriented embedding, introduced in \cite {BLP1}, definition 2.6.1 :

\begin{defn}

Let $(A,\tau)$ be an orthogonal involution with $A$ of even degree, and let $\nu : \Delta (E) \to Z(A,\tau)$ be an orientation. An embedding $ \iota : (E, \sigma) \to (A, \tau)$ is called an {\it oriented
embedding} with respect to $\nu$ if there exist $a \in F^{\times}$ and $\alpha \in A^{\times}$  satisfying the following conditions :

\medskip

{\rm (a)} ${\rm Int}(\alpha) : (A,\theta_{a}) \to (A,\tau)$   is an isomorphism of algebras with involution such that ${\rm Int}(\alpha)  \circ  \epsilon = \iota$.

\medskip

{\rm (b)} The induced automorphism  $c(\alpha) :  Z(A,\theta_{a})  \to Z(A,\tau)$ satisfies $$c(\alpha) \circ   u_{a} = \nu.$$

\end{defn}

We say that there exists an oriented embedding of algebras with involution with respect to $\nu$ if there exists $ (\iota, a, \alpha)$ as above. The elements $ (\iota, a, \alpha,\nu)$ are called {\it parameters} of the oriented embedding.

\section{The obstruction groups}\label{obstruction groups}

\medskip
{\bf A general construction}

\medskip

Recall from \cite {B} the following construction. Let $I$ be a finite set, and let $C(I)$ be the set of maps
$I \to {\bf Z}/2{\bf Z}$. Let $\sim$ be an equivalence relation on $I$. We denote by $C_{\sim}(I)$ the set
of maps that are constant on the equivalence classes. Note that $C(I)$ and $C_{\sim}(I)$ are finite
elementary abelian 2-groups. 

\medskip
{\bf An example}

\medskip
This example will be used in \S \ref{indep}. We say that two finite extensions
$K_1$ and $K_2$ of a field $K$ are  {\it independent} if the tensor product $K_1 \otimes_K K_2$ is a field. 
Let $E = \underset{i \in I} \prod E_i$ be a product of finite field extensions $E_i$ of $K$, and let us consider the
equivalence relation generated by

\medskip
\centerline {$i \sim j \iff$ $E_i$ and $E_j$ are independent over $K$.}

\medskip
Let $C_{\rm indep}(E)$ be the quotient of $C_{\sim}(I)$ by the constant maps; this is a finite
elementary abelian 2-group.

\medskip
{\bf Commutative \'etale algebras with involution}

\medskip
Let $(E,\sigma)$ be a commutative \'etale $L$-algebra with involution such as in \S \ref{definitions}. Note that $E$
is a product of fields, some of which are stable by $\sigma$, and the others come in pairs,  exchanged by $\sigma$. Let us
write $E = E' \times E''$, where 
$E' = \underset{i \in I} \prod E_i$ with $E_i$ a field stable by $\sigma$ for all $i \in I$, and where $E''$ is a product
of fields exchanged by $\sigma$. With the notation of \S \ref {definitions}, we have $F = F' \times F''$, with
$F' = \underset{i \in I} \prod F_i$, where $F_i$ is the fixed field of $\sigma$ in $E_i$ for all $i \in I$. 
Note that $E'' = F'' \times F''$.
For all $i \in I$, let $d_i \in F_i^{\times}$ be such that $E_i = F_i (\sqrt d_i)$, and let $d =(d_i)$.

\medskip
{\bf The subsets $V_{i,j}$}

\medskip
Let $V$ be a set, and for all $i,j \in I$ let $V_{i,j}$ be a subset of $V$. We consider the equivalence
relation on $I$ generated by $i \sim j \iff V_{i,j} \not = \varnothing$.

\medskip
{\bf Global fields}

\medskip
Assume that $K$ is a global field, and let $V_K$ be the set of places of $K$. If $v \in V_K$, we denote by $K_v$ the 
completion of $K$ at $v$.

\medskip
For all $i \in I$, let $V_i$ be the set of places $v \in V_K$ such that there exists a place of $F_i$ above $v$ that
is inert or ramified in $E_i$. For all $i,j \in I$, set $V_{i,j} = V_i \cap V_j$, and let $\sim$ be the
equivalence relation generated by $i \sim j \iff V_{i,j} \not = \varnothing$. 

\medskip
Let $C(E,\sigma)$ be the quotient of $C_{\sim}(I)$ by the constant maps; note that $C(E,\sigma)$ is a finite
elementary abelian 2-group.

\medskip
As a consequence of \cite {BLP1}, Theorem 5.2.1, we show the following. Let  $(A,\tau)$ be as in \S \ref{definitions}. For all $v \in V_K$, set $E^v = E \otimes_K K_v$
and $A^v = A \otimes_K K_v$. 

\begin{theo}\label{0}
Assume that for all $v \in V_K$, there
exists an oriented embedding $(E^v,\sigma) \to (A^v,\tau)$, and that $C(E,\sigma) = 0$. Then there
exists an embedding $(E,\sigma) \to (A,\tau)$. 

\end{theo}

The proof is given in section \ref {ns}, as a consequence of 
Theorem \ref {main}.

\section{Embedding data}

Let $K$ be a global field, and let $(E,\sigma)$ and $(A,\tau)$ be as above. The aim of this section is to recall
some notions from \cite {BLP1} that we need in the following section.

\medskip
We start by recalling from \cite {BLP1} the notion of  embedding data. Assume that for all
$v \in V_K$ there exists an  embedding 
$(E^v,\sigma) \to (A^v,\tau)$. The set of $(a) = (a^v)$, with $a^v \in (F^{v})^{\times}$,  such that for all $v \in V_K$  we have
$(A_v,\tau) \simeq (A_v,\theta_{a^v})$, is called a {\it local embedding datum}. This is sufficient if $(A,\tau)$
is unitary or split orthogonal; however, when $(A,\tau)$ is nonsplit orthogonal, we need the notion of {\it oriented local embedding
data}, as follows. 

\medskip
Let us introduce some notation.

\begin{notation} If $M$ is a field, let  ${\rm Br}(M)$ be the Brauer group of $M$. For $a, b \in M^{\times}$, we
denote by $(a,b)$ the class of the quaternion algebra determined by $a$ snd $b$ in ${\rm Br}(M)$.

\medskip
For $K$ and $F$ as above, and for $v \in V_K$, set $F^v = F \otimes _ K K_v$, and we denote by ${\rm cor}_{F^v/K_v}
: {\rm Br}(F^v) \to {\rm Br}(K_v)$ the corestriction map. Recall that we have a homomorphism ${\rm inv}_v : {\rm Br}(K_v)
\to {\bf Q}/{\bf Z}$. 

\end{notation}

\medskip
{\bf Oriented embedding data}

\medskip Assume that $(A,\tau)$ is nonsplit orthogonal, and let $\nu : \Delta (E) \to Z(A,\tau)$ be an orientation. Suppose that for all $v \in V_K$ there exists an 
oriented  embedding 
$(E^v,\sigma) \to (A^v,\tau)$. An {\it oriented embedding datum} will be an element $(a) = (a^v)$ with $a^v \in (F^{v})^{\times}$ such that for all $v \in V_K$ there exists $\alpha^v \in (A^v)^{\times}$  such that 
$({\rm Int}(\alpha) \circ \epsilon, a^v,\alpha^v,\nu)$ are parameters of an oriented embedding, and that moreover
the following conditions are satisfied
see \cite {BLP1}, 4.1) :

\smallskip
\noindent
$\bullet$ Let $V'$ be the set of places $v \in V_K$ such that 
$\Delta (E^v)$ is a quadratic extension of $K_v$. Then  ${\rm cor}_{F^v/K_v}(a^v,d) = 0$
for almost all $v \in V_K''$.

\smallskip
\noindent
$\bullet$ We have $$\underset {v \in V_K} \sum  {\rm cor}_{F^v/K_v}(a^v,d)  = 0.$$

\medskip
We denote by ${\mathcal L}(E,A)$ the set of oriented local embedding data (of course, the orientation
is only required in the nonsplit orthogonal case - if $(A,\tau)$ is unitary or split orthogonal, then 
${\mathcal L}(E,A)$ is by definition the set of local embedding data).

\section{A necessary and sufficient condition}\label{ns}

Let $K$ be a global field, and let $(E,\sigma)$ and $(A,\tau)$ be as in the previous sections. The aim of this section is to reformulate
the necessary and sufficient condition of \cite {BLP1}; the only difference is a simpler description of the obstruction 
group. Suppose that for all $v \in V_K$ there exists an 
oriented  embedding 
$(E^v,\sigma) \to (A^v,\tau)$, and let $(a) = (a^v_i) \in \mathcal L (E,A)$ be an oriented local embedding datum. 
For all $i \in I$, recall that  $d_i \in F_i^{\times}$ is such that $E_i = F_i (\sqrt d_i)$. 

\medskip Let $C(E,\sigma)$ be the
group defined in \S \ref{obstruction groups}. 
We define a homomorphism $\rho = \rho_a : C(E,\sigma) \to {\bf Q}/{\bf Z}$ by setting

$$\rho_a(c) = \underset {v \in V_K} \sum c(i)  \ {\rm inv}_v {\rm cor}_{F^v/K_v} (a_i^v,d_i).$$

We have the following

\begin{theo}\label{main}
{\rm (a)} The homomorphism $\rho$ is independent of the choice of $(a) = (a^v_i) \in \mathcal L (E,A)$.

\medskip
{\rm (b)} There exists an embedding of algebras with involution $(E,\sigma) \to (A,\tau)$ if and only if
$\rho = 0$.

\end{theo}

\noindent
{\bf Proof.} As we will see, the theorem follows from \cite {BLP1}, Theorems 4.4.1 and 4.6.1, and from the fact that
the group $C(E,\sigma)$ above and the group $\sha(E',\sigma)$ of \cite {BLP1} are isomorphic, fact
that we prove here. Note that part (a) of the theorem can also be deduced directly from \cite {B}, Proposition 12.4.

\medskip
Let us recall the definition of  $\sha(E',\sigma)$ from \cite {BLP1}, 5.1 and \S 3. Recall that 
$E' = \underset{i \in I} \prod E_i$ with $E_i$ a field stable by $\sigma$, and
$F' = \underset{i \in I} \prod F_i$, where $F_i$ is the fixed field of $\sigma$ in $E_i$ for all $i \in I$. As in \cite {BLP1},
\S 3, let 
$\Sigma_i$ be  the set of $v \in V_K$ such that all the places of $F_i$ above $v$ split
in $E_i$; in other words,  $\Sigma_i$ is the complement of $V_i$ in $V_K$.
Let $m = |I|$. 
Given an m-tuple $x = (x_1,..., x_m)\in({\bf Z}  /2 {\bf Z})^{m}$, set $$I_0 = I_0(x) = \{ i \ | \ x_i = 0 \},$$ $$I_1 = I_1(x) = \{ i \ | \ x_i = 1 \}.$$ 
Let $S'$ be the set
$$ S' = \{(x_1,...,\ x_m)\in ({\bf Z}  /2 {\bf Z})^{m} \ | 
(\underset{i\in I_0}{\cap}\Sigma_i)\cup(\underset{j\in I_1}{\cap}\Sigma_j)=V_K \} ,$$ and let $S = S' \cup(0,...,0)\cup(1,...,1).$
Componentwise addition induces a group structure on $S$ (see \cite {BLP1}, Lemma 3.1.1).
We denote by $\sha(E',\sigma)$ the quotient of  $S$ by the subgroup generated by $(1,\dots,1)$.

\medskip
We next show that the groups $\sha(E',\sigma)$ and $C(E,\sigma)$ are isomorphic. The first remark
is that with the above notation, we have
$$ S' = \{(x_1,...,x_m)\in ({\bf Z}  /2 {\bf Z})^{m} \ | 
(\underset{i\in I_0}{\cup}V_i)\cap(\underset{j\in I_1}{\cup}V_j)= \varnothing \}.$$ 

Let us consider the map $F : ({\bf Z}  /2 {\bf Z})^{m} \to C(I)$ sending $(x_i)$ to the map $c : I \to {\bf Z}/2{\bf Z}$
defined by $c(i)= x_i$. We have

$$ S' = \{(x_1,...,x_m)\in ({\bf Z}  /2 {\bf Z})^{m} \ | 
(\underset{c(i) = 0}{\cup}V_i)\cap(\underset{c(j) = 1}{\cup}V_j)= \varnothing \}.$$ Note that this shows that
the following two properties are equivalent :

\medskip
(a) $(x_i) \in S$

\medskip
(b) If $i,j \in I$ are such that $V_i \cap V_j \not = \varnothing$, then we have   $c(i) = c(j).$

\medskip
Recall from \S 2 the definition of the group $C(E,\sigma)$. We consider the equivalence
relation on $I$ generated by
$i \sim j \iff V_{i,j} \not = \varnothing$, and we denote by $C_{\sim}(I)$ the set of $c \in C(I)$ that
are constant on the equivalence classes.

\medskip
Since $(a) \implies (b)$, $F$ sends $S$ to $C_{\sim}(I)$. Moreover, $F$ is clearly injective. Let
us show that $F : S \to C_{\sim}(I)$ is surjective : this follows from the implication $(b) \implies (a)$. 
Hence we obtain an isomorphism of groups $S \to C_{\sim}(I)$, inducing an isomorphism of
groups $\sha(E',\sigma) \to C(E,\sigma)$, as claimed. 

\medskip
The isomorphism $F : \sha(E',\sigma) \to C(E,\sigma)$ transforms $\overline f : \sha(E',\sigma) \to {\bf Q}/{\bf Z}$
defined in \cite {BLP1}, 5.1 and 4.4 into the homomorphism $\rho : C(E,\sigma) \to  {\bf Q}/{\bf Z}$ defined above.
By \cite {BLP1}, Theorem 4.4.1 (see also 5.1) the homomorphism $\overline f$ is independent of the
choice of $(a) = (a^v_i) \in \mathcal L (E,A)$, and this implies part (a) of the theorem. Applying
Theorem 4.6.1 and  Proposition 5.1.1, we obtain part (b).

\medskip
\noindent
{\bf Proof of Theorem \ref {0}.} This is an immediate consequence of Theorem \ref {main}.

\section {An application - independent extensions}\label{indep}

We keep the notation of the previous section; in particular, $K$ is a global field.  We say that two finite extensions
$K_1$ and $K_2$ of $K$ are  {\it independent} if the tensor product $K_1 \otimes_K K_2$ is a field [if $K_1$ and $K_2$ are
subfields of a field extension $\Omega$ of $K$, then this means that $K_1$ and $K_2$ are linearly disjoint). The first
result of the section is the following :  

\begin{theo}\label{2} Assume that $E = E_1 \times E_2$, where $E_1$ and $E_2$ are independent field extensions
of $K$, both stable by $\sigma$. Then $C(E,\sigma) = 0$.

\end{theo}

\noindent
{\bf Proof.} Let us show that $V_1 \cap V_2 \not = \varnothing$. Let $\Omega/K$ be a Galois extension containing $E_1$
and $E_2$, and set $G = {\rm Gal}(\Omega/L)$. Let $H_i \subset G_i$ be subgroups of $G$ such that for $i = 1,2$,
we have $E_i = \Omega^{H_i}$ and $F_i = \Omega^{G_i}$. Since $E_i$ is a quadratic extension of $F_i$, the subgroup
$H_i$ is of index 2 in $G_i$. By hypothesis, $E_1$ and $E_2$ are linearly disjoint over $K$, therefore
$[G : H_1 \cap H_2 ] = [G : H_1] [G : H_2]$. Note that $F_1$ and $F_2$ are also linearly disjoint over $K$, hence
$[G : G_1 \cap G_2] = [G : G_1][G : G_2]$.
This implies that $[G_1 \cap G_2: H_1 \cap H_2] = 4$, hence 
the quotient $G_1 \cap G_2/H_1 \cap H_2$ is an elementary abelian group of order 4. 

\medskip
The field $\Omega$ contains the composite fields $F_1F_2$ and $E_1E_2$. By the above argument, the
extension $E_1E_2/F_1F_2$ is biquadratic. Hence there exists a place of $F_1F_2$ that is inert in both
$E_1F_2$ and $E_2F_1$. Therefore there exists a place $v$ of $K$ and places $w_i$ of $F_i$ above $v$ 
that are inert in $E_i$ for $i = 1,2$.

\medskip

Let $(A,\tau)$ be a central simple algebra as in the previous sections.

\begin{coro} Assume that $E = E_1 \times E_2$, where $E_1$ and $E_2$ are independent field extensions
of $K$, both stable by $\sigma$, and suppose that for all $v \in V_K$ there exists an 
oriented  embedding 
$(E^v,\sigma) \to (A^v,\tau)$. Then there exists an embedding of algebras with involution $(E,\sigma) \to (A,\tau)$.

\end{coro}

\noindent
{\bf Proof.} This follows from Theorems \ref{main} and \ref{2}.

\medskip
To deal with the case where $E$ has more than two factors, we introduce a group $C_{\rm indep}(E,\sigma)$. 
As in \S \ref {obstruction groups}, we write $E = E' \times E''$, where 
$E' = \underset{i \in I} \prod E_i$ with $E_i$ a field stable by $\sigma$ for all $i \in I$, and where $E''$ is a product
of fields exchanged by $\sigma$. Let $\approx$ be the equivalence relation on $I$ generated by

\medskip
\centerline {$i \approx j \iff$ $E_i$ and $E_j$ are independent over $K$.}

\medskip
We denote by $C_{\rm indep}(E,\sigma) = C_{\rm indep}(E')$ the group constructed in \S \ref {obstruction groups} using this
equivalence relation.

\begin{theo}\label {disjoint group} Assume that $C_{\rm indep}(E,\sigma) = 0$, and
suppose that for all $v \in V_K$ there exists an 
oriented  embedding 
$(E^v,\sigma) \to (A^v,\tau)$. Then there exists an embedding of algebras with involution $(E,\sigma) \to (A,\tau)$.

\end{theo}

\noindent
{\bf Proof.} Recall that the group $C(E,\sigma)$ is constructed using the equivalence relation
generated by $i \sim j \iff V_i \cap V_j \not = \varnothing$. Theorem \ref {2} implies that if $E_i$ and $E_j$ are independent over $K$, then $V_i \cap V_j \not = 
\varnothing$, hence $i \equiv j \implies i \sim j$. By hypothesis, $C_{\rm indep}(E,\sigma) = 0$, therefore
$C(E,\sigma) = 0$; hence Theorem \ref{0} implies the desired result. 

\begin{coro}\label{1} Assume that there exists $i \in I$ such that $E_i$ and $E_j$ are independent
over $K$ for all $j \in I$, $j \not = i$. Suppose that for all $v \in V_K$ there exists an 
oriented  embedding 
$(E^v,\sigma) \to (A^v,\tau)$. Then there exists an embedding of algebras with involution $(E,\sigma) \to (A,\tau)$.

\end{coro}

\noindent
{\bf Proof.} Since there exists $i \in I$ such that $E_i$ and $E_j$ are independent
over $K$ for all $j \in I$, $j \not = i$, the group  $C_{\rm indep}(E,\sigma)$ is trivial,  hence the result
follows from Theorem \ref{disjoint group}.

\begin{coro}\label{pairwise} Assume that the fields $E_i$ are pairwise independent over $K$.
Suppose that for all $v \in V_K$ there exists an 
oriented  embedding 
$(E^v,\sigma) \to (A^v,\tau)$. Then there exists an embedding of algebras with involution $(E,\sigma) \to (A,\tau)$.

\end{coro}

\noindent
{\bf Proof.} This follows immediately from Corollary \ref{1}.

\bigskip
\bigskip
Eva Bayer--Fluckiger 

EPFL-FSB-MATH

Station 8

1015 Lausanne, Switzerland

\medskip

eva.bayer@epfl.ch

\end{document}